\documentclass[a4paper,12pt]{article}
\usepackage{amssymb}
\usepackage{amsmath}
\usepackage{graphicx}
\usepackage{epstopdf}

\newcommand{\re}[1]{(\ref{#1})}
\def\mp{{{\cal M}^+_{\lambda,\Lambda}}}
\def\mm{{{\cal M}^-_{\lambda,\Lambda}}}

\newcommand{\eneqa}{\end{eqnarray}}
\newcommand{\begeqaet}{\begin{eqnarray*}}
\newcommand{\eneqaet}{\end{eqnarray*}}
\newcommand{\be}{\begin{equation}}
\newcommand{\ee}{\end{equation}}

\newcommand{\rn}{\rbig^N}

\newcommand{\rbig}{{\mathbb{R}}}

\newcommand{\opr}{{\Omega^\prime}}
\newcommand{\oprb}{{\overline{\Omega^\prime}}}

\newcommand{\ed}{\end{document}}

\newcommand{\ep}{\varepsilon}

\def\beq{\begin{equation}}
\def\eeq{\end{equation}}

\def\Box{\hfill\framebox(0.25,0.25){}}

\newtheorem{thm}{Theorem}

\newcommand{\begeqa}{\begin{eqnarray}}

\begin{document}

\begin{center}{\bf\Large A new method of proving a priori bounds \\ for superlinear elliptic PDE}
\end{center}\smallskip

 \begin{center}
Boyan SIRAKOV\footnote{e-mail : bsirakov@mat.puc-rio.fr}\\
PUC-Rio, Departamento de Matematica,\\ Gavea, Rio de Janeiro - CEP 22451-900, BRAZIL
\end{center}\bigskip

{\small \noindent{\bf Abstract}. We describe a new method of proving a priori bounds for positive supersolutions and solutions of superlinear elliptic PDE, based on global weak Harnack  inequalities and a quantitative Hopf lemma. Novel results based on the method include: (i) equations without a boundary condition on the whole boundary; (ii) equations with nonlinearities which do not have precise growth at infinity; (iii) systems of inequalities with opposite sign. }
%
%\noindent 2000 Mathematics Subject Classification : 35J60, 35J65,
%35P30.
% }

\section{Introduction and Main Results}

In this work we exhibit a new method of proving uniform estimates for positive solutions of nonlinear elliptic PDE, in divergence or non-divergence form, with superlinear growth in the unknown function.

A fundamental and very extensively studied problem in the theory of elliptic PDE is whether nonnegative solutions of equations like
\begin{equation}\label{princ}
-F[u] = f(x,u)
\end{equation}
are uniformly bounded provided the behaviour of $f$ with respect to $u$ for large or small values of $u$ is different from that of the elliptic operator $F$. In this work $F$ will be supposed to be positively $1$-homogeneous in $u$, that is, to have a linear growth in $u$. Then the problem is particularly challenging when $f$ grows more rapidly, i.e. when $f$ is superlinear in $u$  at infinity.

Here is a brief summary of the novelties in our results.
\begin{itemize}
\item[(i)] We do not assume that $f$ has precise power growth as $u\to\infty$, and we do not assume that the equation is satisfied in a domain with a Dirichlet condition on the whole boundary; that is, we obtain both local and global bounds. To our knowledge, all other results in the literature contain at least one of these two hypotheses. Specifically for non-divergence form equations, all previous results assume that $f$ behaves like pure power as $u\to\infty$.
\item[(ii)] We obtain separate results for sub- and super-solutions which can be combined into a $L^\infty$-bound. Thus, we can consider not only equations, but also systems of inequalities with opposite sign, satisfied by different functions, and containing different operators and nonlinearities.
\item[(iii)] The proofs are strikingly short.
\end{itemize}

The method consists in consecutive use of three fundamental estimates in regularity theory: the quantitative maximum principle, the weak Harnack inequality and the local maximum principle. In order to consider equations in a domain, we need the recent up-to-the-boundary extensions of these estimates from  \cite{Si1} and \cite{Si2}. It is worth observing that  $L^\infty$-bounds are traditionally viewed as being independent of regularity theory; one usually proves an a priori bound in $L^\infty$ and then uses it by referring to regularity estimates, in order to infer an uniform bound in higher-order Sobolev or H\"older spaces. We see that there are situations when $L^\infty$-bounds and regularity can be developed starting from the same fundamentals.
A comprehensive discussion of our method, its advantages and shortcomings, will be given in Section \ref{disc}.\smallskip

\smallskip

 We next give precise statements. We assume we have a domain $\Omega\subset\mathbb{R}^n$, $n\ge2$, and functions $b,h\ge0$ which belong to some Lebesgue space $L^q(\Omega)$, with $q>n$. We denote with $T$  a $C^{1,1}$-smooth relatively open portion of the boundary of $\Omega$,  and let $\opr$ be a bounded domain such that $\oprb\subset\Omega\cup T$. We set $B=\|b\|_{L^q(\Omega)}$, $H=\|h\|_{L^q(\Omega)}$, and denote $d(x) = \mathrm{dist}(x,T)$. We also fix constants $\Lambda\ge \lambda>0$.

In what follows $f(x,s)$, $g(x,s)$ will  be  nonnegative functions, measurable in $x\in \Omega$ and continuous in $s\ge0$.
We will assume that (superlinearity)
\begin{equation}\label{superlin}
\lim_{s\to\infty}\frac{f(x,s)}{s}=\infty
\end{equation}
uniformly in $x\in \omega$ for some  subset $\omega\subset \opr$ with positive measure; and that for some $a_0,r>0$,
\begin{equation}\label{hypog}
 {g(x,s)}\le {a_0}{(1+s^r)},  \qquad x\in\Omega,\;s\ge0.
 \end{equation}

Let $\xi:[0,\infty)\to[0,\infty)$ be a increasing bijection such that
\begin{equation}\label{hgrow}
\limsup_{s\to\infty} \frac{\xi(s)}{s^{\beta}}<\infty,
\end{equation}
for some $\beta>0$. Many applications would need only $\xi(s)=s$, $\beta=1$.

We introduce two constants which play an important role. We denote with $p^*$ the exact upper bound of the set of all $p>0$ such that if $u$ is a positive supersolution of $-Lu\ge h$ in $\Omega$ then $u^p\in L^1(\opr)$, for any  second-order operator
$Lu=a_{ij}(x)\partial _{ij}u + b_i(x)\partial_i u $
where $A(x)=(a_{ij}(x))$ is a matrix with eigenvalues in $[\lambda,\Lambda]$ and $\|b_i\|_{L^q(\Omega)}\le B$. We denote with $n_*$ the exact lower bound  of all $p>0$ such that if $c,h\in L^p(\Omega)$ and $u$ is a subsolution of $-Lu-cu\le h$ in $\Omega$, $u\le0$ on $T$, then $u$ is locally bounded in $\opr\cup T$ in terms of the norms of the coeffcients.

In other words, $p^*$ is the optimal exponent in the global weak Harnack inequality, while $n_*$ is the optimal integrability of zero order terms for the validity of the local maximum principle (or the ABP inequality). We know that $p^*>0$ exists and that $n_*\in [n/2,n]$. Explicit expressions for these constants are not known for general operators (and this is an important open problem in the theory of equations in non-divergence form); however if the principal part of $L$ can be written in divergence form as in \re{divsystineq1}-\re{divsystineq2} below, we know that $p^*=n/(n-1)$ and $n_*=n/2$. If in addition $\opr\subset\subset \Omega$ then $p^*=n/(n-2)_+$. See Section \ref{prel} for details.

\begin{thm}\label{apreq}
Let $u\in C(\Omega\cup T)$ be a nonnegative viscosity solution of
\begin{eqnarray}
-\mm(D^2u) +b(x)|Du|&\ge& f(x,u) - h(x) \label{gensystineq1} \\
-\mp(D^2(\xi(u)))-b(x)|D(\xi(u))| &\le& g(x,\xi(u))+ h(x), \label{gensystineq2}
\end{eqnarray}
in $\Omega$, with $u=0$ on $T$. Here $f$ satisfies \re{superlin},  and  $g$ satisfies \re{hypog} for some
\begin{equation}\label{subcrit}
r<1+ \frac{p^*}{\beta n_*}.
\end{equation}
Then for some  $C$ depending on $n, \lambda, \Lambda, B, q, r, |\omega|, \mathrm{diam}(\opr), \|T\|_{C^{1,1}}, f, H,\xi$, and dist$(\opr,\partial\Omega\setminus T)$ if $\opr\not=\Omega$,
 $$
u(x)\le C,\qquad\xi(u(x))\le C d(x), \qquad x\in\Omega^\prime.
$$
%where $C$ depends on $n, \lambda, \Lambda,  B, a_0, s_0, \opr, \mathrm{dist}(\opr,\partial\Omega\setminus T),f, H$.
\end{thm}

In the left hand side of \re{gensystineq1}--\re{gensystineq2} we have Pucci's extremal operators, i.e. the supremum and infimum of linear second order operators with fixed bounds for the coefficients.  They can of course be replaced by arbitrary (and different in the two inequalities) linear operators in the form
$-a_{ij}(x)\partial _{ij}u + b_i(x)\partial_i u $
where $A(x)=(a_{ij}(x))$ has eigenvalues in $[\lambda,\Lambda]$, $\|b_i\|_{L^q(\Omega)}\le B$.

As we already mentioned, $p^*$ and $n_*$ can be specified for divergence form operators (such as the Laplacian). Because of the importance of this case, we state the result separately.

\begin{thm}\label{aprdiveq} Let $A^{(i)}=\left(a^{(i)}_{ij}(x)\right)$, $i=1,2$, be matrices such that $A^{(i)}\ge \lambda I$, $\|A^{(i)}\|_{W^{1,q}(\Omega)}\le \Lambda$.
 Assume that $u, \xi(u)\in H^1(\Omega)$ satisfy in the weak sense
\begin{eqnarray}
-\mathrm{div}\left( A^{(1)}(x)Du \right) + b(x) |Du|&\ge& f(x,u) - h(x)\label{divsystineq1}\\
-\mathrm{div}\left( A^{(2)}(x)D(\xi(u)) \right) - b(x) |D(\xi(u)) | &\le&  g(x,\xi(u))+ h(x),\label{divsystineq2}
\end{eqnarray}
and $u=0$ on $T$. Here $f$ satisfies \re{superlin},  and  $g$ satisfies \re{hypog} for some
\begin{equation} r<\frac{n+1}{n-1}+\left(\frac{1}{\beta}-1\right)\frac{2}{n-1},
 \label{divsubcrit}
\end{equation}
or $\displaystyle r<\frac{n}{(n-2)_+}+\left(\frac{1}{\beta}-1\right)
 \frac{2}{(n-2)_+}$ if  $\opr\subset\subset \Omega$.
 Then for some $C$ as above
 $$
u(x)\le C,\qquad\xi(u(x))\le C d(x), \qquad x\in\Omega^\prime.
$$

\end{thm}

In \cite{So} Souplet constructed an example which shows that positive solutions of the Dirichlet problem for $-\Delta u=a(x) u^{r}$ in a smooth bounded domain are not uniformly bounded, if $r>\frac{n+1}{n-1}$, $0\lneqq a(x)\le a_0$. This implies that  the upper bound in \re{divsubcrit} cannot be improved. \smallskip

\noindent {\it Remark 1}. In each of the above two theorems we consider the appropriate notion of a weak solution --  in the viscosity sense (see \cite{CIL}, \cite{CCKS}), or in the weak integral sense (see \cite[Chapter 8]{GT}). In the setting of Theorem \ref{aprdiveq} the two notions are essentially equivalent, see \cite{Si2} and its appendix.\smallskip

\noindent {\it Remark 2}.  In the case $\opr\subset\subset \Omega$ Theorem \ref{aprdiveq} is valid for any uniformly elliptic $A^{(i)}\in L^\infty(\Omega)$. However, only boundedness of $A^{(i)}$ is not sufficient for the proof of the  bounds up to the boundary, since even basic boundary behavior results such as Hopf lemma fail for uniformly elliptic operators in divergence form with only bounded (or even continuous) coefficients.  \smallskip

\noindent {\it Remark 3}. Larger ranges for the exponent $r$ are known when considering {\it classical solutions of specific equations}. For instance, there is an a priori bound   for classical solutions of the Dirichlet problem for $-\Delta u= u^{r}$ (or equations which blow up to this equation at each point, see \cite{GS}) if $r<\frac{n+2}{n-2}$, due in particular to the conformal invariance of the Laplacian.\smallskip

\noindent {\it Remark 4}. We can allow unbounded dependence in $x$; for instance $a_0$ in \re{hypog} can be divided by an appropriate power of $d(x)$, or $a_0$ can be a function in some $L^q(\Omega)$ for sufficiently large $q$. This amounts to using the H\"older inequality in the proofs, and is left to the interested reader. \smallskip

\noindent {\it Remark 5}. A model situation in which we need $\xi(u)\not= u$ is when we start with an equation for which we seek an a priori bound, and notice there are two changes of the unknown function, one of which is a subsolution, the other a supersolution, to inequalities as in  \re{gensystineq1}-\re{gensystineq2} or \re{divsystineq1}-\re{divsystineq2}.\smallskip

Our method separates estimates for super- and sub-solutions. In particular, we obtain an uniform integrability estimate for supersolutions, worth to be stated independently.

\begin{thm}\label{aprineq}
Assume $u\in C(\Omega)$ is a nonnegative viscosity supersolution  of
\begin{equation}\label{gensupineq}
-\mm(D^2u) + b(x)|Du| \ge f(x,u) - h(x)\qquad\mbox{in } \;\Omega.
\end{equation}
Then for some $p^*$ and $\ep^*$ depending
 on $n, \lambda, \Lambda,B,$ and any $p<p^*$ and $\ep<\ep^*$,
\begin{equation}\label{resaprineq}
\left(\int_{\opr} u^p\right)^{1/p}\le C,\qquad \left(\int_{\opr}
\left(\frac{u}{d}\right)^{\varepsilon}\right)^{1/\varepsilon}\le C,
\end{equation}
with $C$  depending on $n, \lambda, \Lambda, B, p, q, \ep, |\omega|, \mathrm{diam}(\opr), \|T\|_{C^{1,1}}, f, H$, as well as dist$(\opr,\partial\Omega\setminus T)$ if $\opr\not=\Omega$. 

If in the left-hand side of \re{gensupineq} we have an operator in divergence form as in Theorem \re{aprdiveq}, the estimate \re{resaprineq} holds for every nonnegative supersolution $u\in H^1_{\mathrm{loc}}(\Omega)$, and all $\varepsilon<1$ and $\displaystyle p<\frac{n}{n-1}$ (or $\displaystyle p<\frac{n}{n-2}$ if $\opr\subset\subset \Omega)$.
\end{thm}

We stress that Theorem \ref{aprineq} is valid for supersolutions, without a requirement for $u$ to be a   subsolution (let alone a solution) of an equation, and that no boundary condition for $u$ is needed.  The only requirement on $f$ is the superlinearity condition \re{superlin}; in particular ``indefinite" nonlinearities are allowed, i.e. $f$ may vanish on a nontrivial subset of the domain $\Omega$.\smallskip

An instance of the method described here, in a somewhat more complicated form and restricted  to \re{gensystineq1} with logarithmically superlinear nonlinearities, was  sketched in the note \cite{Si3} and then used in \cite{NS}.

\section{On methods for proving a priori bounds}\label{disc}

The importance of obtaining an uniform estimate for positive solutions to a nonlinear elliptic problem was recognized in the early stages of the development of the theory, and fundamental contributions were obtained in the nowadays classical papers by Brezis and Turner \cite{BT}, Gidas and Sprick \cite{GS}, de Figueiredo, Lions and Nussbaum \cite{dFLN}. We also refer to the survey \cite{L1}.

The available methods of proving a priori bounds for positive solutions of elliptic boundary value problems can be divided, broadly, into two categories. First, methods that require and use variational structure of the equation, i.e. its equivalence to an integral equality involving test functions; and second, methods which require and use scaling invariance of the equation, most often  the fact that $u$ and $t^{-a}u(tx)$, $t>0$, satisfy the same type of equation for some $a>0$, at least asymptotically as $t\to \infty$.

The results in the first group are global in the sense that they  need an equation set in a bounded domain with a Dirichlet condition prescribed on the whole boundary. For instance, this is required in an important step in any use of these methods, namely, testing the equation with the first eigenfunction of the operator, which leads to a bound in $L^1_d$ for the solution $u$ and the nonlinearity $f(x,u)$ ($L^1_d$ is the space of functions whose product with the distance to the boundary is in $L^1$).

For an equation to have variational structure it is necessary that the elliptic operator be in divergence form. Methods based on use of such structure include:
\begin{itemize}
\item The method of Brezis-Turner \cite{BT}, in which the $L^1_d$ estimate for $f(x,u)$ together with the growth assumption for $f$ and Hardy-Sobolev inequalities lead to a $H^1$ estimate which is bootstrapped into an $L^\infty$ one.
\item  The method of de Figueiredo-Lions-Nussbaum \cite{dFLN}. The $L^1_d$ estimate for $f(x,u)$ together with the growth assumption for $f$ and Pohozaev inequalities lead to a $H^1$ estimate which is bootstrapped into an $L^\infty$ one, through a Moser type trick which bounds the $H^1$-norm of a power of $u$.  Similarly to the methods based on scaling, better assumptions on the growth of $f$ can be achieved, but at the price of requiring either convexity of the domain  or a global hypothesis on the behaviour of $f$ (to make possible using moving planes).
\item The method of Quittner-Souplet \cite{QS}, which uses weighted Lebesgue spaces. The $L^1_d$ estimate for $f(x,u)$ together with the growth assumption for $f$ and $L^{p_1}_d$-to-$L^{p_2}_d$ inequalities for the equation, $p_1<p_2$,  lead to a  bootstrap for the $L^p_d$-norm of $u$, ending into an $L^\infty$ bound. This method represents a nontrivial and important extension of the one in~\cite{BT},  as it allows for weaker (``very weak") solutions and leads to much better results when applied to systems of equations.
    \end{itemize}

Our method improves the above in that we obtain local estimates, i.e. we do not need a full domain with a boundary condition. However, we cannot deal with very weak solutions, that is, functions which are only in $L^1_d$ rather than in the Sobolev space $H^1$, since the results in the next section are not available for such solutions (this is an interesting open problem in itself).\medskip

The scaling, or ``blow-up", method was introduced by Gidas and Spruck in \cite{GS}. It requires a hypothesis of specific behaviour of $f$ as $u\to\infty$, such as $\lim_{s\to\infty} f(\lambda s)/f(s) = \lambda^p$ for some $p>0$ and all $\lambda >0$. The scaling method proceeds always by contradiction, if an a priori bound is not available, a rescaling of an exploding sequence of solutions at or close to points of maximum converges to a positive solution of a simplified equation in the whole or a half-space. The conclusion is then obtained through a non-existence theorem for that simplified equation. The passage to the limit requires a regularity estimate (at least a $C^\alpha$-estimate) for the elliptic operator.

An important extension of the scaling method of \cite{GS} was obtained by  Polacik-Quittner-Souplet \cite{PQS}. They use a topological doubling lemma which dispenses the need of a boundary condition, leading to an uniform estimate for the solutions in terms of a (negative) power of the distance to the boundary. The rescaling is done around points provided by the doubling lemma around which the function cannot double its  size.

In the scaling method, most frequently the simplified limit equation is $-\Delta u = u^p$ which, by the results in \cite{GS2}, does not admit entire positive solutions for $p<(n+2)/(n-2)$ (as opposed to $(n+1)/(n-1)$). This is why when the scaling method is applicable, it leads to a priori bounds for a larger range of growth of the nonlinearities. However, the nonexistence result of~\cite{GS2} hinges both on the conformal invariance of the Laplacian and the specific form of the nonlinearity $u^p$ (for developments in this direction we refer to~\cite{LN} and the references there). Since our method is based solely on properties shared by all uniformly elliptic operators, and is independent of the form of the nonlinearity, it does not appear to be amenable to obtain results for nonlinearities with power growth above $(n+1)/(n-1)$. We recall the latter is the threshold for a priori bounds for a general equation, as shown in \cite{So}.

To summarize, when reduced to solutions of a scalar equation, our method improves the results coming from  variational methods in that we do not need a boundary condition on the whole boundary; and the results coming from scaling methods in that we do not need a precise behaviour of $f$ for large values of $u$. To our knowledge, our result here is the first whatsoever for non-divergence form operators, in which precise power growth on $f$ is not assumed. On the negative side, our method does not allow to obtain stronger results in terms of the growth of $f$ in $u$, in the particular cases when such results are available thanks to the specific nature of the operator, the domain, and the nonlinearity. We do not have results for very weak solutions, in so far as the results in the next section are not proved for such solutions.\smallskip

Another essential feature of our method, which appears to be new, is that we treat separately subsolutions and supersolutions. Combining these separate results gives a $L^\infty$ bound. Thus we can prove a priori bounds for two different inequalities of opposite sign, satisfied by different functions depending on each other; for instance, we have almost free of charge the possibility of considering arbitrary $\xi$, resp. $\beta>0$, in the main theorems.\smallskip

Furthermore, our argument is direct, not by contradiction. No passage to the limit is involved, nor any use of regularity estimates is required.\smallskip

We note that we allow indefinite nonlinearities, that is, $f(\cdot,u)$ is allowed to vanish in some subdomain. Previous results with indefinite nonlinearities are essentially restricted to problems with scaling invariance,  see for instance \cite{BCN}, \cite{AL}, \cite{So}, \cite{GIQ}, \cite{R},  \cite{DL}. These works contain various nontrivial extensions of the scaling method;  many of them allow also functions in $x$ that change sign, a situation which we do not consider here. Of course,  our results  can actually treat some sign-changing nonlinearities, since depending on the hypotheses, we can apply our bounds separately to the positive or negative part of $u$, in the subdomains where the nonlinearity has a sign.

\section{Preliminaries: global estimates for subsolutions and supersolutions}\label{prel}

In the proof of our results we use the following global estimates from \cite{Si1}, \cite{Si2}. They are up-to-the-boundary extensions of the classical interior bounds by de Giorgi-Moser (for equations in divergence form) and  Krylov-Safonov (for equations in non-divergence form). For a discussion on their ramifications, use, and history, we refer to \cite{Si1}-\cite{Si2}, where a large set of references is available.

 In this section we assume we have a bounded $C^{1,1}$-domain $G\subset\rn$, and set $d(x)=\mathrm{dist}(x,\partial G)$. The constants $c,C$ in the theorems below can depend on the diameter of $G$ and an upper bound for the curvature of $\partial G$.

The following theorem is a quantification of the Hopf lemma.

\begin{thm}\label{bqsmp} (\cite[Theorem 1.1]{Si1}) Assume that $u$ is a nonnegative viscosity supersolution of
 \begin{equation}\label{ineq0}
-\mm(D^2u) + b(x)|Du| \ge  h(x) \quad\mbox{in } \;G.
\end{equation}
Then there exist constants $\varepsilon,c,C>0$ depending on $n$, $\lambda$, $\Lambda$, $ q$,  and $\|b\|_{L^q(G)}$, such that
 \begin{equation}\label{ineq1}
\inf_{G} \frac{u}{d} \ge c\left( \int_{G} (h^+)^{\varepsilon}\right)^{1/{\varepsilon}} - C\|h^-\|_{L^q(G)}.
\end{equation}
\end{thm}

The next result is a global extension of the  weak Harnack inequality.

\begin{thm}\label{bwhi} (\cite[Theorem 1.2]{Si1}, \cite[Theorem 1.1]{Si2}, \cite[Corollary 1.1]{Si2}) Assume that $u$ is a nonnegative viscosity supersolution of \re{ineq0}. Then there exist constants $p^*, \varepsilon^*>0$ depending on $n$, $\lambda$, $\Lambda$,  $q$ and $\|b\|_{L^q(G)}$, such that for each $\ep<\ep^*$, $p<p^*$,
 \begin{equation}\label{ineq2}
\max\left\{ \left( \int_{G} \left(\frac{u}{d}\right)^\varepsilon\right)^{1/\varepsilon},\left(\int_{G} u^p\right)^{1/p}\right\} \le C\left( \inf_{G} \frac{u}{d} + \|h^+\|_{L^q(G)}\right).
\end{equation}
If $u$ is a supersolution of
\begin{equation}\label{ineq00}
-\mathrm{div}(A(x)Du) +b(x)|Du|\ge   h(x)\quad\mbox{in } \;G,
\end{equation}
then
$$p^*=\frac{n}{n-1}\,,\qquad \ep^*=1.
$$
\end{thm}

We also record the following Lipschitz bound, which is a boundary extension of the local maximum principle for subsolutions. We write $x=(x^\prime, x_n)\in \mathbb{R}^{n-1}\times \mathbb{R}$ and denote the cube  $Q_R= Q_R^\prime \times (0,R)$, where $Q_R^\prime=\{ -R/2<|x^\prime|<R/2\}\subset \mathbb{R}^{n-1}$. We also set $Q_R^0 = Q_R^\prime \times \{0\}$, the lower boundary of $Q_R$.

\begin{thm}\label{blmp}
 There exists $n_*\in [n/2, n)$ such that if $p>n_*$, $c,h\in L^p(Q_2)$,  $r>0$, $u$ is a subsolution of
  \begin{equation}\label{ineq000}
-\mp(D^2u) -b(x)|Du| -c(x)u\le   h(x) \quad\mbox{in } \;Q_2,\qquad
u\le0\quad\mbox{on } \;Q_2^0,
\end{equation}
there is $C>0$ depending on  $n$, $\lambda$, $\Lambda$, $p$, $q$, $r$,  $\|b\|_{L^q(Q_2)}$, $\|c\|_{L^p(Q_2)}$, such that
 \begin{equation}\label{ineq320}
\sup_{Q_1} {u^+}\le C\left(\left( \int_{Q_2} (u^+)^r\right)^{1/r} + \|h^+\|_{L^p(Q_2)}\right).
\end{equation}
If in addition $c,h\in L^q(Q_2\cap\{x_n<\delta\})$ for some $\delta>0$, $q>n$, then for some constant $C$ depending also on the $L^q$-norm of $c$  and $\delta$,
 \begin{equation}\label{ineq32}
\sup_{Q_1} \frac{u^+}{x_n}\le C\left(\left( \int_{Q_2} (u^+)^r\right)^{1/r} + \|h^+\|_{L^p(Q_2)} + \|h^+\|_{L^q(Q_2\cap\{x_n<\delta\})}\right).
\end{equation}

If instead $u$ is a subsolution of $$-\mathrm{div}(A(x)Du) -b(x)|Du| -c(x)u\le  h \quad\mbox{in } \;Q_2,\qquad
u\le0\quad\mbox{on } \;Q_2^0,
$$
 then the above holds with $\displaystyle n_*=\frac{n}{2}$.
\end{thm}

\noindent {\it Proof.} The interior version of this result can be found in \cite[Theorem 8.17]{GT}, \cite[Theorem 3.1]{KS3}. Extending it up to the boundary is not difficult, since $u^+$ extended as $0$ outside $Q_2$ is a subsolution in $Q_2^\prime\times(-2,2)$, and the interior estimate can be applied. Combining this with $C^1$-bounds gives \re{ineq2}; the full argument can be found in the proof of \cite[Theorem 1.3]{Si1} which is \re{ineq32} with $n_*=n$ and $c=0$. In order to get the same for  $c\not=0$ it is enough to consider the term $cu$ as a right-hand side, i.e. to replace $h$ by $h+cu$, apply \re{ineq32} with $q$ replaced by $q-\delta$, use the H\"older inequality to write
$$
\|cu\|_{L^{q-\delta}}\le \|c\|_{L^q}\|u\|_{L^{A_\delta}}
$$
for sufficiently small $\delta>0$ (so that $q-\delta>n_*$) and $A_\delta\in \mathbb{R}$, and finally downgrade the $L^A$-norm of $u$ to a $L^r$-norm through a well-known analysis argument on a shrinking sequence of cubes (given for instance on pages 75-76 of \cite{HL}).

\section{Proofs of the a priori bounds}

We recall we write $x=(x^\prime, x_n)\in \mathbb{R}^{n-1}\times \mathbb{R}$, and  $Q_R= Q_R^\prime \times (0,R)$, where $Q_R^\prime=\{ -R/2<|x^\prime|<R/2\}\subset \mathbb{R}^{n-1}$.

By charting and locally straightening $T$, as well as scaling, we can assume that  $ Q_4\subset\Omega$ and $Q_4^\prime\times\{0\}\subset T$. By a covering argument, it is enough to prove the estimates in $Q_1$.

Note that the change of variables which straightens the boundary changes a solution to any of \re{gensystineq1}, \re{gensystineq2}, \re{divsystineq1}, \re{divsystineq2} into a solution of a similar inequality, with possibly modified $\lambda, \Lambda, B$, depending only on the $C^{1,1}$-norm of $T$.

\subsection{Proof of Theorem \ref{aprineq}}

  Let $G$ be a smooth domain, with the $C^2$-norm of its boundary being a universal constant, such that $Q_2\subset G\subset Q_4$; and $d(x)=\mathrm{dist}(x,\partial G)$. Fix $k=|\omega|/2>0$ and $\omega^\prime = \omega\cap\{x_n>k\}$ (so that $|\omega^\prime|\ge |\omega|/2=k$).

Set
$$
A_u=\inf_{x\in Q_1} \frac{u(x^\prime,x_n)}{x_n}.$$
 By applying  Theorem \ref{bqsmp} to the inequality \re{gensupineq} we obtain
 $$
 A_u\ge \inf_{x\in G} \frac{u(x)}{d(x)}\ge c\left(\int_G f^\ep(x,u(x))\,dx\right)^{1/\ep}-C\|h\|_{L^q(G)},
 $$
hence
\begeqaet
(A_u+ C\|h\|_{L^q(Q_4)})^\ep&\ge&c \int_{Q_1} f^\ep(x,u(x^\prime,x_n))\,dx_n\,dx^\prime \\
&=& c\int_{Q_1^\prime}\int_0^1 f^\ep(x,x_n\frac{u(x^\prime,x_n)}{x_n})\,dx_n\,dx^\prime \\
&\ge&c \int_{Q_1^\prime}\int_0^1 \inf_{s\in[A_u,\infty)} f^\ep(x,sx_n) \, dx_n\,dx^\prime \\
&\ge& c\int_{Q_1^\prime}\int_{k}^1 \inf_{s\in[A_u,\infty)} f^\ep(x,sx_n) \, dx_n\,dx^\prime \\
&\ge& c\int_{\omega^\prime} \inf_{s\in[kA_u,\infty)} f^\ep(x,s) \, dx \\
&\ge& c|\omega^\prime|\left(\inf_{x\in \omega^\prime,s\in[kA_u,\infty)} f(x,s)\right)^\ep.
\eneqaet
Thus, recalling $|\omega^\prime|\ge k$,
$$
 \inf_{x\in \omega,s\in[kA_u,\infty)} f(x,s) \le   (c k^{1+\ep} )^{-1/\ep} (kA_u) + C\|h\|_{L^q(Q_4)} =C\left( kA_u + \|h\|_{L^q(Q_4)}\right),
$$
which implies $kA_u\le C$, by the superlinearity assumption \eqref{superlin}. That is, $A_u\le C$,
where $C$ depends on the right quantities and is independent of $u$.

Theorem \ref{aprineq} then follows from  Theorem \ref{bwhi} applied to  $G$, since
$$
\left(\int_{Q_1} \left(\frac{u}{x_n}\right)^\varepsilon\right)^{1/\varepsilon}\le \left(\int_{G} \left(\frac{u}{d}\right)^\varepsilon\right)^{1/\varepsilon} \le C A_u+ C\|h\|_{L^q(Q_4)}\le C.\qquad \Box
$$

\subsection{Proof of Theorem \ref{apreq}}

Let   $C, s_0>0, s_1\ge s_0$, be such that
$$
g(x,s)\le 2a_0 s^r\;\mbox{ if }\;s\ge s_0,\qquad\mbox{and}\qquad s_0\le \xi(s) \le Cs^\beta\;\mbox{ if }\; s\ge s_1.
$$
From \re{gensystineq2} it is easy to see that the function
$$
\tilde u = \max\{u, s_1\}
$$
is a weak (viscosity) solution of
$$
-\mp(D^2(\xi(\tilde u)))-b(x)|D(\xi(\tilde u))| \le g(x,\xi(\tilde u))+ h(x)
$$
since the maximum of subsolutions is a subsolution.

We set
$$
0\le c(x) := \frac{ g(x,\xi(\tilde u(x))) }{ \xi(\tilde u(x))}\le 2 a_0\xi(\tilde u)^{r-1}\le C\tilde u^{\beta(r-1)}.
$$
We now use the hypothesis $n_*\beta(r-1)<p^*$  and Theorem \ref{aprineq}, in which we showed that $\tilde u\in L^p(Q_2)$ for each $p<p^*$, to deduce that $c(x)$ is uniformly bounded in $L^{p_0}(Q_2)$, for some $p_0>n_*$.

We set $v=\xi(\tilde u)-\xi(s_1)$ and observe that
$$
-\mp(D^2v)-b(x)|Dv| - c(x) v \le h(x) + \xi(s_1) c(x),
$$
in $Q_2$, so by Theorem \ref{blmp}, inequality \re{ineq320} (applied with $Q_{3/2}$ instead of $Q_1$ and with some fixed $r<p^*/\beta$ so that $\xi(\tilde u)$ be uniformly bounded in  $L^r(Q_2)$), we can infer that $v$ is bounded in $L^\infty(Q_{3/2})$. The latter is clearly equivalent to $\xi(u)$ being uniformly bounded in $L^\infty(Q_{3/2})$.

Hence the function $\tilde h(x)= g(x,\xi(u))+ h(x)$ is in $L^q(Q_{3/2})$, $q>n$, and by Theorem \ref{blmp}, inequality \re{ineq32}, with $Q_{3/2}$ instead of $Q_2$, applied to $$-\mp(D^2(\xi( u)))-b(x)|D(\xi(u))| \le \tilde h(x),$$
we conclude
$$
\left\|\frac{\xi(u)}{x_n}\right\|_{L^\infty(Q_1)}\le C.
$$

This proves Theorem \ref{apreq}. To prove Theorem \ref{aprdiveq} it is enough to observe that \re{divsubcrit} is \re{subcrit} with $p^*=n/(n-1)$ and $n_*=n/2$ (see also Theorems \ref{bwhi}-\ref{blmp}). The upper bound for the exponent in the case $\opr\subset\subset\Omega$ comes from replacing $n/(n-1)$   by $n/(n-2)_+$, since the latter is the optimal $p^*$ in the {\it interior} weak Harnack inequality, \cite[Theorem 8.18]{GT}, and we can cover $\opr$ with balls in the interior of $\Omega$.

\end{document}